\documentclass{elsart}
\usepackage{amsmath,amssymb} 
\usepackage[all]{xy}

\newcommand{\N}{\mbox{$\mathbb N$}}
\newcommand{\Q}{\mbox{$\mathbb Q$}}
\newcommand{\PTree}{\mbox{${\rm PTree}$}}
\newcommand{\PTreeM}{\mbox{${\rm PTree}\{(M_k)_{k\in\N}\}$}}

\newcommand{\PRTrx}{\mbox{${\rm PRTree}\{M^X\}$}}

\newcommand{\PRTree}{\mbox{${\rm PRTree}$}}
\newcommand{\YTree}{\mbox{${\rm YTree}$}}
\newcommand{\emp}{\mbox{$\emptyset$}}
\newcommand{\sh}{\mbox{${\sqcup\!\sqcup}$}}
\newcommand{\As}{\mbox{${\mathcal As}$}}
\newcommand{\Com}{\mbox{${\mathcal Com}$}}
\newcommand{\Cmg}{\mbox{${\mathcal Cmg}$}}
\newcommand{\Mag}{\mbox{${\mathcal Mag}$}}
\newcommand{\Mgom}{\mbox{${\mathcal Mag}_{\omega}$}}
\newcommand{\Mgtwo}{\mbox{${\mathcal Mag}_{2}$}}
\newcommand{\Mgn}{\mbox{${\mathcal Mag}_{N}$}}
\newcommand{\Mgnprime}{\mbox{${\mathcal Mag}_{N'}$}}
\newcommand{\Cmgn}{\mbox{${\mathcal Cmg}_{N}$}}
\newcommand{\Dend}{\mbox{${\mathcal Dend}$}}
\newcommand{\Lie}{\mbox{${\mathcal Lie}$}}
\newcommand{\Sab}{\mbox{${\mathcal Sab}$}}
\newcommand{\Brace}{\mbox{${\mathcal Brace}$}}
\newcommand{\Prim}{\mbox{${\rm Prim}$}}
\DeclareMathOperator{\id}{id}

\def\Y{\setlength{\unitlength}{.4pt}\begin{picture}(60,40)(0,0)
\put(30,0){\line(0,1){10}} \put(30,10){\line(-1,1){30}}
\put(30,10){\line(1,1){30}}
\end{picture}}

\begin{document}
\begin{frontmatter}

\title{On Hopf algebra structures over free operads}

\author{Ralf Holtkamp}
\address{University of Bochum, 44780 Bochum, Germany}
\ead{ralf.holtkamp@ruhr-uni-bochum.de}

\date{}

\begin{abstract}
The operad Lie can be constructed as the operad of primitives
 Prim$\As$ from the operad $\As$ of associative algebras. This is
 reflected by the theorems of Friedrichs, Poincar\'e-Birkhoff-Witt and
 Cartier-Milnor-Moore. We replace the operad $\As$ by families of
 free operads ${\mathcal P}$, which include the operad $\Mag$ freely generated by a noncommutative
 non-associative binary operation and the operad of Stasheff
 polytopes. We obtain Poincar\'e-Birkhoff-Witt type theorems and
 collect information about the operads Prim${\mathcal P}$, e.g.\
 in terms of characteristic functions.
\end{abstract}
\begin{keyword}
Operad \sep Hopf algebra \sep primitive \sep
Poincar\'e-Birkhoff-Witt

\MSC 16W30 \sep 17A50 \sep 18D35 \sep 18D50
\end{keyword}
\end{frontmatter}

\section*{Introduction}

Recent developments in Hopf algebra theory show that there are
important objects which could be called "non-classical" Hopf
algebras.
 Typical examples are
dendriform Hopf algebras \cite{lilrhopftree,liroa,lirobcr,lirob}.
Dendriform algebras (introduced in \cite{lilodia}) are equipped
with two operations $\prec$, $\succ$ whose sum is an associative
multiplication. Many of these Hopf algebras provide new links
between algebra, geometry, combinatorics, and theoretical physics,
e.g.\ renormalization in quantum field theory in the work of
Connes and Kreimer (cf.\ \cite{lick}). Other examples are magma
Hopf algebras and infinitesimal Hopf algebras (cf.\
\cite{liasinf}). Here the usual Hopf algebra axioms have to be
changed.

 Let $\Mag$ be the operad freely generated by a
non-commutative non-associative binary operation $\vee^2(x_1,x_2)$
also denoted by $x_1\cdot x_2$. Let similarly $\Mgom$ be freely
generated by $n$-ary operations $\vee^n$, one for each $2\leq n\in
\N$. A basis for the space of $n$-ary operations is given by
reduced planar rooted trees with $n$ leaves. This is the operad of
Stasheff polytopes, see \cite{listaphys}. More generally, we
consider operads $\Mgn$, $2\leq N \leq \omega$, intermediating
between $\Mag$ and $\Mgom$.

In a natural way, Hopf algebra structures over these operads can
be considered. In fact it is sufficient that the operad $\mathcal
P$ fulfills certain coherence conditions. Then the free $\mathcal
P$-algebra generated by a set $X$ of variables is always a
$\mathcal P$-Hopf algebra with the diagonal $\Delta_a$ (also
called co-addition) as a comultiplication. There canonically
exists an operad $\Prim{\mathcal P}$, and $\Prim\As=\Lie$.  A
table of pairs ${\mathcal P}, \Prim{\mathcal P}$ (triples in fact,
if coassociativity of the comultiplication is replaced by a
different law) is given by Loday and Ronco in
\cite{lilr03,lilr04}, and the family given by $\Mgn$ might be
added to the list. In the case of $\Mag$, a full set of
multilinear primitives was given by Shestakov and Umirbaev in
\cite{lisu}, answering a problem posed by Hofmann and Strambach in
\cite{lihofstr}. An approach to study the primitive elements via
Taylor expansions with respect to the so-called algebra of
constants was given in \cite{lighmag}.

In order to describe the operads $\Prim\Mgn$, we describe the
graded duals of the given $\Mgn$-Hopf algebras. These duals are
equipped with commutative multiplications $\sh$. We show that
these commutative algebras are freely generated by the primitive
elements (and dually the given generalized Hopf algebras are
connected co-free). This is an analogon of the
Poincar\'e-Birkhoff-Witt theorem.

We conclude that the characteristic of the $\Sigma_n$-module
$\Prim\Mgn(n)$ is of the form
\begin{equation*}
{\rm ch}_n\bigl(\Prim\Mgn(n)\bigr)= \frac{1}{n}\sum_{d\vert
n}\mu(d) c[N]'_{\frac{n}{d}}\ p_d^{\frac{n}{d}}.
\end{equation*}
Especially this means that the dimension of ${\Prim\Mgn}(n)$ is
equal to
\begin{equation*}
(n-1)!\ c[N]'_n.
\end{equation*}
 For $N=2$, $c[2]'_k=c'_k$ is the $k$-th
$\log$-Catalan number, with $c'_1=1,\ c'_2=1,\ c'_3=4,\ c'_4=13,\
c'_5=46,\ c'_6=166,\ \ldots\ $

\medskip

In Section 1, we recall basic facts about trees and operads. We
introduce the notion of admissibly labeled (planar or abstract)
rooted trees, where a sequence of sets $M_k$ contains the allowed
labels for vertices of arity $k$. This is useful to define the
operads $\Mgn$. Related integer sequences and their logarithmic
derivatives are also considered.
\newline
In Section 2, we use the notion of unit actions on operads (see
\cite{lilosci}, and cf.\ \cite{liefguo} for further studies) in a
generalized setting. The properties of the unit actions on $\Mgn$
are needed in Section 3 to define the corresponding ${\mathcal
P}$-Hopf algebras. In their definition we do not include objects
with various not necessarily associative operations and not
necessarily coassociative cooperations. The given definition is
general enough to include dendriform Hopf algebras, though, and we
sketch the relations to work of Loday, Ronco, and others.
\newline
In Section 4 we describe the free unitary $\Mgom$-algebra
$K\{X\}_{\omega}$ together with $\Delta_a$ in more detail. We also
explicitly describe the graded dual $(K\{X\}_{\omega},\Delta_a)$,
which is equipped with a shuffle multiplication that is a sum of
the shuffles. For $\Mgn$-Hopf algebras one may pass to the
appropriate subalgebras and quotients.
\newline
An analogon of the Poincar\'e-Birkhoff-Witt theorem for the
operads $\Mgn$ is proved in Section 5, see Theorems\ \ref{thmpbw}
and\ \ref{thmpbwdual}. Our approach makes use of co-${\bf D}$
object structures, where ${\bf D}$ is the category of unitary
magmas. We also discuss the existence of a Cartier-Milnor-Moore
theorem and of Eulerian idempotents. Moreover, we note that there
are cocommutative $\As$-Hopf algebras (with the same coalgebra
structure) associated to the given $\Mgn$-Hopf algebras.
\newline
In Section 6, we discuss the generating series and characteristic
functions for the operads of primitives. Here we focus on
$\Prim\Mag$ and $\Prim\Mgom$. For small $n$, we can present the
corresponding $\Sigma_n$-modules in terms of irreducible
representations.
\newline
 We also relate our results to recent independent
work of Bremner - Hentzel - Peresi \cite{libhp} and of
P\'erez-Izquierdo \cite{liper} concerning the case of
$\Mag$-algebras. We discuss the description of $\Prim\Mag(4)$ by
Sabinin operations.

This paper is based on my Habilitationsschrift \cite{lihohab}. I
would like to thank F.\ Chapoton, L.\ Gerritzen, J.-L.\ Loday,
J.-C.\ Novelli, and I.\ P.\ Shestakov for helpful discussions and
remarks. I also thank the referee for helpful suggestions.

\begin{section}{Trees and free operads}\label{sectone}

\begin{subsection}{Some combinatorics of trees}

 A finite connected graph $\emptyset\neq T=({\rm
Ve}(T)$, ${\rm Ed}(T))$, with a distinguished vertex
$\rho_T\in{\rm Ve}(T)$, is called an abstract rooted tree (with
root $\rho_T$), if for every vertex $\lambda\in{\rm Ve}(T)$ there
is exactly one path connecting $\lambda$ and $\rho_T$. At each
vertex there are incoming
 edges and exactly one outgoing edge. (Here
we think of the edges as being oriented towards the root, and
 we add to the root
an outgoing edge that is not connected to any further vertex.)
 At a given vertex $\lambda$,
the number $n$ of incoming edges is called the arity
ar$_{\lambda}$ of $\lambda$. We write the set ${\rm Ve}(T)$ of
vertices as a disjoint union $\bigcup_{n\in\N} {\rm Ve}^n(T)$. The
vertices of arity $0$ are called leaves, and we denote ${\rm
Ve}^0(T)$ by ${\rm Le}(T)$.

 An abstract rooted tree $T$
together with a chosen order of incoming edges at each vertex is
called a planar rooted tree (or ordered rooted tree), cf.\
\cite{listan}, \cite{listanb} as a general reference.

 It is well-known that the number of planar rooted trees
with $n$ vertices is the $n$-th Catalan number
$ c_n=\frac{(2(n-1))!}{n!(n-1)!}=\sum_{l=1}^{n-1}c_lc_{n-l}.$
The sequence of Catalan numbers is
$ c_1=1, c_2=1, c_3=2,c_4=5, c_5=14, c_6=42, c_7=132,
 c_8=429, c_9=1430, \ldots $
with generating series $f(t)=\sum_{n=1}^{\infty}c_nt^n$ given by
$\frac{1-\sqrt{1-4t}}{2}.$

The numbers $c_n$ also count the number of planar binary rooted
trees with $n$ leaves (or $2n-1$ vertices). A tree $T$ is called
binary, if ${\rm Ve}(T)={\rm Ve}^2(T)\cup {\rm Le}(T)$.

\medskip

Let $a_n, n\geq 1,$ be a sequence of integers with generating
series $f(t)=\sum_{n=1}^{\infty}a_n t^n\!\!.$ The logarithmic
derivative of $f(t)$ is the series
$g(t):=\frac{\partial}{\partial t}\log\bigl(1+f(t)\bigr),$
and we say that the sequence $a'_n, n\geq 1,$ with
$\sum_{n=1}^{\infty} a'_n t^n=t\cdot g(t)$ is obtained from $a_n,
n\geq 1,$ by logarithmic derivation.

The sequence of $\log$-Catalan numbers $c'_n$, starting with
1, 1, 4, 13, 46, 166, 610, 2269, 8518, 32206, $\ldots$
has the generating series $\frac{2t}{3\sqrt{1-4t}-1+4t}.$

Since
$\frac{2}{3\sqrt{1-4t}-1+4t}=\frac{\partial}{\partial
t}\log\bigl(\frac{3-\sqrt{1-4t}}{2}\bigr),$
it is obtained by logarithmic derivation from the Catalan numbers
$c_n, n\geq 1$.

In the set of all planar (rooted) trees with $n$ vertices,
 the number of vertices with even arity is given
by $c'_n$. The corresponding numbers of vertices with odd arity
have the generating series
\begin{equation*}
\begin{split}
\sum_{n=1}^{\infty} & n c_n t^n - \sum_{n=1}^{\infty} c'_n t^n =
\frac{t}{\sqrt{1-4t}}-\frac{2t}{(3-\sqrt{1-4t})\sqrt{1-4t}}
\\
 &=\frac{t(1-\sqrt{1-4t})}{(3-\sqrt{1-4t})\sqrt{1-4t}}
 =t^2 + 2t^3 + 7t^4 + 24t^5 + 86 t^6 + 314 t^7 + \ldots\\
\end{split}
\end{equation*}
(see \cite{lisl} A026641, \cite{lids} p.\ 258).

\goodbreak

For example, in the set \newline $
 
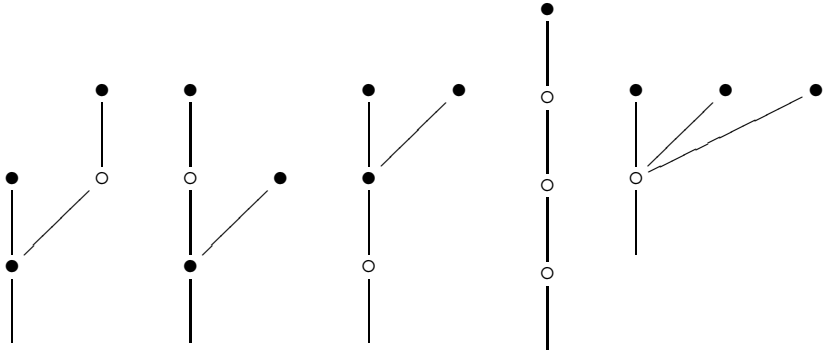
\xymatrix{ \\ & {\bullet}\ar@{-}[d] \\ {\bullet}\ar@{-}[d]& {\circ}\ar@{-}[dl]\\
{\bullet}\ar@{-}[d]\\ \\}\ \ \ \ \ \
 \xymatrix{ \\ {\bullet}\ar@{-}[d] \\ {\circ}\ar@{-}[d]& {\bullet}\ar@{-}[dl]\\
{\bullet}\ar@{-}[d]\\ \\}\ \ \ \ \ \
 \xymatrix{ \\ {\bullet}\ar@{-}[d] & {\bullet}\ar@{-}[dl]\\ {\bullet}\ar@{-}[d]\\
{\circ}\ar@{-}[d]\\ \\}\ \ \ \ \ \
 \xymatrix{{\bullet}\ar@{-}[d] \\ { \circ}\ar@{-}[d] \\
{\circ}\ar@{-}[d]\\ {\circ}\ar@{-}[d]\\
\\}\ \ \ \ \ \
 \xymatrix{\\ {
\bullet}\ar@{-}[d] & { \bullet}\ar@{-}[dl] &
 { \bullet}\ar@{-}[dll]\\ {\circ}\ar@{-}[d]\\ \\}$

of planar trees with $n=4$ vertices, we count seven vertices with
odd arity and $c'_4=13$ vertices with even arity.

\end{subsection}
\begin{subsection}{Admissible Labelings}

Let $M$ be a set and $T$ a planar (or abstract) tree. Then a
labeling of $T$ is a map $\nu:{\rm Ve}(T)\to M$. The tree $T$
together with such a labeling is called a labeled tree.

Let a collection $M_0, M_1, M_2,\ldots$ of
 sets be given, and let $M=\bigcup_{k\in\N}
M_k$. A labeling  $\nu:{\rm Ve}(T)\to M$ of a planar (or abstract)
tree $T$ is called admissible, if the restrictions $\nu\vert{\rm
Ve^k}(T)$ are maps ${\rm Ve^k}(T)\to M_k$, i.e.\ it holds that:
\begin{equation*}
\nu(\lambda)\in M_k \text{ if \ ar}_{\lambda}=k.
\end{equation*}

If $M_1=\emptyset$, only reduced trees can be provided with an
admissible labeling. A tree $T$ is called reduced, if
ar$_{\lambda}\neq 1$ for all $\lambda\in {\rm Ve}(T)$.

The set  of planar rooted trees $T\in\PTree$ with admissible
labeling from $(M_k)_{k\in\N}$ is denoted by $\PTreeM$.

Non-labeled trees occur as trivially labeled trees, i.e.\ in the
case where all $M_k$, $k\in\N$, are given by a one-element set
$\{\circ\}$. Similarly (setting $M_1=\emptyset$), we consider
non-labeled planar reduced trees.

The number $C_n=\#\PRTree^n$ of planar reduced trees with $n$
leaves is called the $n$-th
 super-Catalan number, also called the $n$-th little Schroeder
 number.
The generating series for the super-Catalan numbers is
$\frac{1}{4}(1+t-\sqrt{1-6t+t^2})$ (cf.\ \cite{lisl} A001003). The
first 10 super-Catalan numbers are
1, 1, 3, 11, 45, 197, 903, 4279, 20793, 103049.

Moreover, for fixed $N\in \N$, we consider the sets $\PRTree[N]^n$
of planar reduced trees which have $n$ leaves and fulfill the
property, that for every vertex $\lambda$ the arity ar$_{\lambda}$
is $\leq N$. The corresponding integer sequence is denoted by
$c[N]_n$. Clearly for $N=2$, we get binary trees and the sequence
$c_n$, while for $N\to \infty$, we exhaust all planar reduced
trees. More exactly, $c[N]_n=C_n$ whenever $n\leq N$.

Similar to the definition of log-Catalan numbers, we define the
sequences $C'_n$ and $c[N]'_n$ by logarithmic derivation. One can
check that these sequences are in fact integer sequences. (This
will also follow later on.)

\end{subsection}

\begin{subsection}{Free Operads}

Let $K$ be a field. Let a collection $M=(M_k)_{k\geq 2}$ of sets
be given, and set
 $M_0:=\{\circ\}, M_1:=\emptyset$.
Then the free non-$\Sigma$-operad $\underline{\Gamma(M)}$
generated by the collection $(M_k)_{k\geq 2}$ can be conveniently
described using admissible labelings:

Let $\underline{\Gamma(M)}(0)=0$,
$\underline{\Gamma(M)}(1)=K\cdot\vert$, where $\vert$ is the tree
consisting of the root. The elements of $\underline{\Gamma(M)}(n)$
are all linear combinations of admissibly labeled (necessarily
reduced) planar trees with $n$ leaves.

The operad structure on the sequence of vector spaces
$\underline{\Gamma(M)}(n)$ is determined by $K$-linear
$\circ_{n,i}$-operations
 \begin{equation*}
\circ_{n,i} : \underline{\Gamma(M)}(n)\otimes
\underline{\Gamma(M)}(m)\to \underline{\Gamma(M)}(m+n-1) , \text{
all } n, m\geq 1,\ 1\leq i\leq n.
\end{equation*}
and the unit $\vert$. Here, if $T^1$ is a tree with $n$ leaves,
and $T^2$ is a further tree having $m$ leaves, $T^1 \circ_{n,i}
T^2$ (for $i=1,\ldots,n$) is given by the substitution of $T^2$ in
$T^1$ at the $i$-th leaf (obtained by replacing the specified leaf
of $T^1$ by the root of $T^2$).

\medskip
 Let
$M_2$ consist of one generator $\alpha$, and let $M_k=\emptyset\
(k\geq 3)$. Then all elements of $\underline{\Gamma(M)}$ are
linear combinations of planar binary trees. For $n\geq 1$, we can
identify a basis of $\underline{\Gamma(M)}(n)$ with the set of
(non-labeled) planar binary trees with $n$ leaves. Especially,
$\dim\underline{\Gamma(M)}(n)=c_n$.

The tree $\Y$ corresponds to the binary operation $\alpha$, and we
get ternary operations $\alpha\circ_{2,1} \alpha$ and
$\alpha\circ_{2,2} \alpha$ as compositions.

This is the non-$\Sigma$-operad $\underline{\Mag}$ of
(non-unitary) magma algebras.

\end{subsection}

\begin{subsection}{Grafting operations}

A word $T^1.T^2\ldots T^k$ (or an ordered tuple $(T^1,T^2,\ldots
T^k)$, not necessarily non-empty) of planar trees is called a
planar forest. Given a forest $T^1.T^2\ldots T^k$ of $k\geq 0$
trees, together with a label $\rho\in M_k$, there is a tree
$T=\vee_{\rho}(T^1.T^2\ldots T^k)$ defined by introducing a new
root of arity $k$ and grafting the trees $T^1,\ldots,T^k$ onto
this new root. The new root gets the label $\rho$, and the
specified order determines the order of incoming edges at
$\rho_T$. The tree $T$ is called the grafting of $T^1.T^2\ldots
T^k$ over $\rho$.

The non-$\Sigma$-operad $\underline{\mathcal K}$ of Stasheff
polytopes is the free non-$\Sigma$-operad $\underline{\Gamma(M)}$
generated by a collection $(M_k)_{k\geq 2}$ of one-element sets.
We denote the generator of arity $k$ by $\vee^k$. It corresponds
to the grafting operation $\vee$ restricted to planar forests
consisting of $k$ trees (and their $K$-linear combinations). The
tree symbolizing this operation is called the $k$-corolla. The
cells of the Stasheff polytope (or associahedron) in dimension
$n-2$ can be identified with reduced planar trees with $n$ leaves
and
 form a basis of $\underline{\mathcal K}(n)$.

\end{subsection}

\begin{subsection}{A family of free operads}

For each $N\in\N$, we can also consider the non-$\Sigma$-operad
$\underline\Mgn$ generated by a collection $(M_k)_{k\geq 2}$ given
by $\{\vee^k\}$ for $k\leq N$ and $\emptyset$ for $k>N$. The trees
symbolizing operations are elements of $\PRTree[N]^n$, i.e.\ they
have at most arity $N$ vertices. Clearly
$\underline\Mgtwo=\underline\Mag$. There are obvious inclusion
maps $\underline\Mgn\to \underline\Mgnprime$, $N\leq N'$, defined.
All these are sub-operads of $\underline{\mathcal K}$, which we
will also denote by $\underline\Mgom$.

Let $\Mgn$ be the operads given by the symmetrizations of the
non-$\Sigma$-operads $\underline{\Mgn}$.

 The symmetrization of
a non-$\Sigma$-operad $\underline{\mathcal P}$ is defined by
${\mathcal P}(n)= \underline{\mathcal P}(n)\otimes_K K\Sigma_n$
(all $n$), where $\Sigma_n$ is the symmetric group. The
composition maps are induced by the maps of $\underline{\mathcal
P}$ (and the maps of the operad $\As$ with $\As(n)=K\Sigma_n$).
Thus, e.g., the $\Sigma_n$-module $\Mag(n)$ is given by $c_n$
copies of the regular representation $K\Sigma_n$, and $\Mgom(n)$
is given by $C_n$ copies of $K\Sigma_n$.
Operads that occur as symmetrizations of non-$\Sigma$-operads are
also called regular operads.

\end{subsection}
\end{section}

\begin{section}{Free $\Mgn$-algebras}\label{secttwo}

Let $K$ be a field of characteristic 0, and let
$X=\{x_1,x_2,\ldots\}$ be a finite or countable set of variables.
We consider algebras over the operads $\Mgn$ defined in the
previous section. If not specified, $N$ is allowed to be in
$\N_{\geq 2}\cup\{\omega\}.$

To ensure that the tensor product of ${\mathcal P}$-algebras is
provided with the structure of a ${\mathcal P}$-algebra there are
several approaches possible (cf.\cite{limoe}, \cite{lilosci},
\cite{lipseudo}). The approach of \cite{lilosci} for binary
quadratic operads is generalized in the following. For the purpose
of this paper, we mostly deal with regular operads.

\begin{defn}

Let ${\mathcal P}$ be an operad, ${\mathcal P}(0)=0, {\mathcal
P}(1)=K=K\id$.
\newline
Let a 0-ary element $\eta$ be adjoined to the $\Sigma$-space
${\mathcal P}$ by ${\mathcal P}'(i):=\begin{cases} {\mathcal P}(i)
&: i\geq 1\\ K \eta &: i=0.\\
\end{cases}$
\newline
A unit action on ${\mathcal P}$ is a partial extension of the
operad composition onto ${\mathcal P}'$ in the sense that
 composition maps (fulfilling the
associativity, unitary, and invariance conditions)
$\mu_{n;m_1,\ldots,m_n}$
 are defined on
\begin{equation*}
{\mathcal P}'(n)\otimes {\mathcal P}'(m_1)\otimes \ldots \otimes
{\mathcal P}'(m_n)\to {\mathcal P}'(m) ,\ m:=m_1+\ldots+m_n,
\end{equation*}
for all $m_j\geq 0$ ($j=1,\ldots,n$), for $n\geq 2,\ m>0$ (or
$n\leq 1, m\geq 0$).
\newline
 Given a ${\mathcal P}$-algebra
$\overline{A}$, we can define structure maps for $A:=K 1\oplus
\overline{A}$ such that $\eta\in {\mathcal P}'(0)$ is mapped to
$1\in A$, and $A$ is called a unitary ${\mathcal P}$-algebra.
\end{defn}
\begin{defn}
Let ${\mathcal P}$ be a regular operad, ${\mathcal P}(n)=
\underline{\mathcal P}(n)\otimes_K K\Sigma_n$, where
$\underline{\mathcal P}$ is generated by fixed sets $(M_k)_{2\leq
k\leq N}.$ Let ${\mathcal P}$ be equipped with a unit action.
\newline
 If  there are operations
 $\star_n\in\underline{\mathcal P}(n)
 =\underline{\mathcal P}(n)\otimes_K K\cdot\id\subseteq\mathcal P(n),
 n\leq N$, fulfilling
\begin{equation*}
\begin{split}
\star_n \circ_{n,i}(\eta) &=\star_{n-1}, \text{ all } i\\
 \star_2\circ_{2,1}(\eta) &=\star_2\circ_{2,2}(\eta)=\id,\\
\end{split}
\end{equation*}
then we say that the unit action respects the operations
$\star_n$, and we define
\begin{equation*}
\mu_{n;0,\ldots,0}(\star_n\otimes\eta\otimes\ldots\otimes\eta)=\eta.
\end{equation*}
Such a unit action is called coherent, if for all ${\mathcal
P}$-algebras $\overline{A},\overline{B}$ it holds that
$(\overline{A}\otimes K 1)\oplus (K 1\otimes \overline{B})\oplus
(\overline{A}\otimes\overline{B})$
  is again a ${\mathcal P}$-algebra
 with:
\newline
For all $n$, all $p\in M_n$, all $a_i\in A, b_i\in B$,
\begin{equation*}
\begin{split}
&p(a_1\otimes b_1,a_2\otimes b_2,\ldots,a_n\otimes b_n):=
\star_n(a_1,a_2,\ldots,a_n)\otimes p(b_1,b_2,\ldots,b_n) \\
&\text{ in case that at least one } b_j\in \overline B,\\
 &p(a_1\otimes 1,a_2\otimes
1,\ldots,a_n\otimes 1):= p(a_1,a_2,\ldots,a_n)\otimes 1, \text{ if
} p(a_1,a_2,\ldots,a_n) \\ &\text{ is defined},
\\
\end{split}
\end{equation*}
and the operations on $(\overline{A}\otimes K 1)\oplus (K 1\otimes
\overline{B})\oplus (\overline{A}\otimes\overline{B})$ generated
by these $p$ via compositions (and the $\Sigma$-action) are
obtained by applying the same compositions (or permutations) on
both tensor components correspondingly.

In case an operad ${\mathcal Q}$ is a quotient of a regular operad
${\mathcal P}$ with a coherent unit action, we say that ${\mathcal
Q}$ is equipped with the induced coherent action if the analogous
equations hold for the images of the operations $\star_n$ and $p$
in ${\mathcal Q}$.
\end{defn}

Since $\Mgn$ is a regular operad freely generated by operations
$\vee^k$, $k\leq N$, we get:

\begin{lem}

 Each operad $\Mgn$ is equipped with a (unique) unit action which
respects the operations $\vee^k$, i.e.\
\begin{equation*}
\begin{split}
\vee^k \circ_{k,i}(\eta) &=\vee^{k-1}, \text{ all } i\\
 \vee^2\circ_{2,1}(\eta) &=\vee^2\circ_{2,2}(\eta)=\id.\\
\end{split}
\end{equation*}
This  unit action is coherent.
\qed
\end{lem}

\begin{rem}{\rm

Let $M^X=(M^X_k)_{k\geq 0}$ denote the sequence
\begin{equation*}
M^X_0=X,\ \  M^X_k=M_k=\{\vee^k\} \text{ for } 2\leq  k\leq N, \ \
\text{ and } M^X_k=\emp \text{ else.}
\end{equation*}
Then the free $\Mgn$-algebra has the set of admissibly labeled
planar trees as a vector space basis. In the case $N=\omega$ these
are the reduced planar trees with leaves labeled by $X$, and for
$N\in\N$ only trees with maximum arity $N$ occur.
\newline
The free $\Mgn$-algebra is naturally graded, such that the planar
trees with $n$ leaves are homogeneous of degree $n$. We set
$K\{X\}_{N}^{(0)}= K 1$ and identify $1$ with the empty tree
$\emp$ (which we now adjoin to the set of reduced planar trees).
For $N\in\N_{\geq 2}\cup\{\omega\}$, we denote the free
$\Mgn$-algebra with unit $1$ by $K\{X\}_{N} =
{\bigoplus_{n=0}^{\infty}} K\{X\}_{N}^{(n)}.$ The admissibly
labeled trees with only one vertex are identified with $X$ and
form a basis of $K\{X\}_{N}^{(1)}$. We also set
$\vee^k(1,1,\ldots,1)=1$ (all $k$).
\newline
Consequently, the free $\Mag$-algebra with unit $1$ can be
identified with the space of labeled binary trees $K\{X\}\subset
K\{X\}_{N}$ (all $N\geq 2$), equipped with the free binary
operation $\cdot=\vee^2$. On $\Mag$, the coherent unit action is
just the one given by $1\cdot a= a\cdot 1 = a$ (for every
$\Mag$-algebra $A$, every $a\in A$), with component-wise
multiplication on the tensor product. The operads $\As$ and $\Com$
are equipped with the induced coherent unit actions.

}
\end{rem}

\end{section}

\begin{section}{${\mathcal P}$-Hopf algebras and $\Prim{\mathcal P}$}

Since Hopf algebras combine operations and cooperations, there is
no operad whose algebras or coalgebras are Hopf algebras. To
describe them, and also generalizations with not necessarily
associative operations and not necessarily associative
cooperations, one would use PROPs. Here we are only interested in
the case where the set of cooperations is generated by one
coassociative cooperation. Thus we do not need the generality of
PROPs and stay close to operad theory.

\begin{defn}

Let ${\mathcal P}$ be an operad equipped with a coherent unit
action. Let $A=K\ 1 \oplus\overline{A}$ be a unitary ${\mathcal
P}$-algebra and $A\otimes A= K\ 1 \oplus (\overline{A}\otimes K
1)\oplus (K 1\otimes \overline{A})\oplus
(\overline{A}\otimes\overline{A})$ be equipped with its unitary
${\mathcal P}$-algebra structure (see Section \ref{secttwo}).
\newline
 Let $\Delta: A\to A\otimes A$
 be a $K$-linear coassociative map (called comultiplication map),
such that $\Delta(1)=1\otimes 1$ and
$\Delta'(a):=\Delta(a)-a\otimes 1 - 1\otimes a \in {\overline
A}\otimes {\overline A}$ for all $a\in {\overline A}$.
\newline
 Then $A$ together with $\Delta$ is called an (augmented) ${\mathcal
P}$-bialgebra, if $\Delta$ is a morphism of unitary ${\mathcal
P}$-algebras, i.e.\ if $\Delta\circ p_A=p_{A\otimes
A}\circ(\Delta\otimes\ldots\otimes\Delta),$ for all $p\in
{\mathcal P}$.

Let moreover $A$ be given by $\cup_{n\in \N} A_n$, where each
$A_i$ is a subspace of $A_{i+1}$. Then $A$ is called a filtered
${\mathcal P}$-Hopf algebra if $\Delta(A_n)\subseteq
\sum_{i=0}^{n} A_i\otimes A_{n-i}  \text{ (all } n).$

We call $A$ connected graded, if $A=\bigoplus_{n\in \N}A^{(n)}$
for finite dimensional vector spaces $A^{(i)}\subseteq A_i$ such
that $A^{(0)}=K$ and
\begin{equation*}
\Delta'(A^{(n)})\subseteq \sum_{i=1}^{n-1} A^{(i)}\otimes
A^{(n-i)} \text{ (all } n\geq 1).
\end{equation*}

\end{defn}

In case that ${\mathcal P}$ is a regular operad equipped with a
coherent unit action, we are going to define a comultiplication
map on $A=K\oplus F_{\mathcal P}(V_X)$. Let us recall that for any
operad ${\mathcal P}$, any sequence $(g_1,g_2,\ldots)$ of elements
of a ${\mathcal P}$-algebra $A'$ uniquely determines a ${\mathcal
P}$-algebra morphism $\gamma=\gamma_{(g_1,g_2,\ldots)}: A\to A'$
with $x_i\mapsto g_i$ (all $x_i\in X$) and $1\mapsto 1$.

\begin{lem}\label{lemdiagonal}

Let ${\mathcal P}$ be a regular operad equipped with a coherent
unit action and let $A$ be the free unitary ${\mathcal P}$-algebra
generated by the vector space $V_X$ with basis $X$. Then there is
a coassociative ${\mathcal P}$-algebra morphism $\Delta_a$ defined
by
\begin{equation*}
\Delta_a(x)=x\otimes 1 + 1\otimes x, \text{ all } x\in X.
\end{equation*}
It provides the free unitary ${\mathcal P}$-algebra $A$ with the
structure of a connected graded ${\mathcal P}$-Hopf algebra.

If we furthermore assume that the $K$-linear map $\tau: A\otimes
A\to A\otimes A, a_1\otimes a_2\mapsto a_2\otimes a_1$ is a
${\mathcal P}$-algebra isomorphism, then $\Delta_a$ is
cocommutative in the sense that $\tau\circ\Delta_a=\Delta_a$.
\end{lem}
\begin{pf}
By construction the map $\Delta_a$ is a ${\mathcal P}$-algebra
morphism  and $\Delta'_a(f)\in {\overline A}\otimes{\overline A}$
for all $f\in {\overline A}$. To check coassociativity, we just
have to verify that
\begin{equation*}
(\Delta_a\otimes\id)\bigl(\Delta_a(x)\bigr)=(\id\otimes\Delta_a)\bigl(\Delta_a\bigr(x)\bigr),
\text{ all } x\in X.
\end{equation*}
Both sides are equal to $x\otimes 1\otimes 1+1\otimes x\otimes
1+1\otimes 1\otimes x$.
\newline
The free unitary ${\mathcal P}$-algebra $A$ allows a grading (with
$A^{(0)}=K$) respected by all ${\mathcal P}$-operations. Similarly
the  ${\mathcal P}$-algebra operations on $A\otimes A$ respect the
grading (component-wise), and the map $\Delta_a$ provides $A$ with
the structure of a connected graded ${\mathcal P}$-Hopf algebra.
\newline
The last assertion follows directly from the definition of
$\Delta_a$.
 \qed
\end{pf}

\begin{rem}\label{remcocom}{\rm

We call $\Delta_a$ the diagonal or co-addition (in analogy to the
 usage in \cite{libh}, where any
abelian group-valued representable functor leads to a co-addition
 on the representing object).

The cocommutative Hopf algebra given by the free unitary
$\As$-algebra $K\langle X\rangle$ together with the diagonal
$\Delta_a$ is well-known.
\newline
 Dually one can equip the standard tensor
coalgebra (i.e.\ $K\langle X\rangle$ with its deconcatenation
coalgebra-structure) with the shuffle multiplication
$\sh=\Delta_a^{*}$, which is commutative (cf.\ \cite{lireu},
Chapter 1).

}
\end{rem}

\begin{defn}

For any connected graded ${\mathcal P}$-Hopf algebra
$A=\bigoplus_{n\in \N} A^{(n)}$, we define the vector space
$A^{*g}=\bigoplus_{n\in \N} (A^{*g})^{(n)}=\bigoplus_{n\in \N}
(A^{(n)})^*$,
where $V^*={\rm Hom}_K (V,K)$ for any vector space $V$. We call
$A^{*g}$ the graded dual of $A$.

We denote by $\Delta^*:A^{*g}\otimes A^{*g}\to A^{*g}$ the
$K$-linear map given by\newline $\bigl(\Delta^*(f_1\otimes
f_2)\bigr)(a)=\sum f_1(a_{(1)})f_2(a_{(2)})$, where $f_1, f_2\in
A^{*g}$ and $a\in A$ with $\Delta(a)=\sum a_{(1)}\otimes a_{(2)}$.
\end{defn}

The maps $\Delta^*$ and $\Delta$ are adjoint with respect to the
canonical bilinear form $\langle, \rangle: A^{*g}\times A\to K$,
i.e.\
\begin{equation*}
\langle \Delta^*(f_1\otimes f_2),a\rangle=\langle f_1\otimes
f_2,\Delta(a)\rangle,
\end{equation*}
where $\langle f_1\otimes f_2, g_1\otimes g_2\rangle = \langle
f_1, g_1\rangle\otimes\langle f_2,g_2\rangle$.

\begin{lem}\label{lemdualhopf}

The space $A^{*g}$ together with the operations $p^*$,
$p\in{\mathcal P}$, is a ${\mathcal P}$-coalgebra.
\newline
The space $A^{*g}$ together with $\Delta^*:A^{*g}\otimes A^{*g}\to
A^{*g}$ is a graded $\As$-algebra with unit $1\in
(A^{*g})^{(0)}=K$.
 Furthermore, a
cocommutative $\Delta$ leads to a commutative multiplication
$\Delta^*$.

\end{lem}

\begin{pf}
This is just the classical result, modified using the fact that a
${\mathcal P}$-algebra structure on $A$ induces a
 ${\mathcal P}$-coalgebra structure on $A^{*g}$.
\qed
\end{pf}

We consider primitive elements, i.e.\ the elements $f$ such that
$\Delta(f)=f\otimes 1 + 1\otimes f$, with respect to the
comultiplication $\Delta_a$ on the free ${\mathcal P}$-algebra
$F_{\mathcal P}(V_X)$ (over an operad ${\mathcal P}$ with coherent
unit action).
\newline
The following Lemma reflects that the composition of primitive
elements is again primitive and gives the definition of the operad
of primitives (compare also \cite{lighmag,lilosci}). We recall
that the space of elements in $F_{\mathcal
P}(V_{\{x_1,x_2,\ldots,x_n\}})^{(n)}$
 which are multilinear (i.e.\ have degree 1 with
respect to each variable $x_i$) can be identified with ${\mathcal
P}(n)$, for any operad ${\mathcal P}$.

\begin{lem}\label{lemprimoperad}

Let ${\mathcal P}$ be as in Lemma \ref{lemdiagonal}. The
$\Sigma_n$-spaces $\Prim{\mathcal P}(n)$ of multilinear primitive
elements in $F_{\mathcal P}(V_{\{x_1,x_2,\ldots,x_n\}})^{(n)}$
define a sub-operad $\Prim{\mathcal P}$ of ${\mathcal P}$, with
free algebra functor $F_{\Prim{\mathcal P}}=\Prim F_{\mathcal P}$.
\end{lem}
\begin{pf}
1) Let $A=K\oplus F_{\mathcal P}(V_X)$ together with $\Delta_a$ be
the connected graded ${\mathcal P}$-Hopf algebra defined in Lemma
\ref{lemdiagonal}. The image of an element $p(x_1,x_2,\ldots)$ of
$A$ under $\Delta_a$ is given by
\begin{equation*}
 \Delta_a\bigl( p(x_1,x_2,\ldots)\bigr)=p_{A\otimes A}
\bigl(x_1\otimes 1 + 1\otimes x_1,\ x_2\otimes 1 + 1\otimes
x_2,\ldots\bigr).
\end{equation*}
This can be expanded into a sum. To simplify notation, we treat
the case where $p\in M_n$ is a generating operation. Two of the
summands are
 $p(x_1,x_2,\ldots)\otimes 1$ and
$\star_n(1,1,\ldots)\otimes p(x_1,x_2,\ldots)=1\otimes
p(x_1,x_2,\ldots)$, because the unit action is coherent. Let
$n\geq 2$, and let us collect all summands with first tensor
component $x_1$, say. Therefore we have to compute
\begin{equation*}
p_{A\otimes A} \bigl(x_1\otimes 1, 1\otimes x_2, 1\otimes x_3,
\ldots,1\otimes
x_n\bigr)=\underbrace{\star_n(x_1,1,\ldots,1)}_{x_1}\otimes
p(1,x_2,x_3,\ldots,x_n).
\end{equation*}
Similarly, in the $\Delta_a$-image of every multilinear
$p(x_1,x_2,\ldots,x_n)$, the term $x_j\otimes
p(\ldots,x_{j-1},1,x_{j+1},\ldots)$ collects all summands with
first tensor component $x_j$. Moreover this implies that every
multilinear primitive $p(x_1,x_2,\ldots,x_n)$ is mapped to 0 if
any $k$ variables are replaced by 1, $1\leq k\leq n-1$.
\newline
2) Given a sequence $(g_1,g_2,\ldots)$ of primitive elements of
$A$, we consider the ${\mathcal P}$-algebra morphism
$\gamma=\gamma_{(g_1,g_2,\ldots)}: A\to A$ with $x_i\mapsto g_i$
(all $x_i\in X$). Then
\begin{equation*}
\bigl(\Delta_a\circ\gamma\bigr)(x_i)=g_i\otimes 1 + 1\otimes
g_i=\bigl((\gamma\otimes\gamma)\circ\Delta_a\bigr)(x_i).
\end{equation*}
 Now let $p(x_1,x_2,\ldots,x_n)$ be a multilinear primitive element
 and let
 $g:=\gamma\bigl(p(x_1,x_2,\ldots,x_n)\bigr)$. We claim that $g$ is
 primitive. Since $\Delta_a$ and $\gamma$ are ${\mathcal P}$-algebra morphisms,
\begin{equation*}
\begin{split}
\Delta_a(g)=&\bigl(\Delta_a\circ\gamma\circ p\bigr)
(x_1,x_2,\ldots)=p_{A\otimes A}
\bigl(\Delta_a(g_1),\Delta_a(g_2),\ldots\bigr)\\ =&p_{A\otimes A}
(g_1\otimes 1 + 1\otimes g_1,\ g_2\otimes 1 + 1\otimes
g_2,\ldots).\\
\end{split}
\end{equation*}
As in step 1), we can expand into a sum where two of the summands
are $\underbrace{p(g_1,g_2,\ldots)}_g\otimes 1$ and
$1\otimes\underbrace{p(g_1,g_2,\ldots)}_g$. The other summands
yield 0 by step 1), because
$\gamma\bigl(p(1,x_2,x_3,\ldots,x_n)\bigr)$,$\ldots$,
$\gamma\bigl(p(1,1,\ldots,1,x_n)\bigr)$ are 0.
\newline
3) The vector space $\Prim{\mathcal P}(n)$ is a $\Sigma_n$-space:
This follows from step 2), when we consider the maps $\gamma$ for
$(g_1,g_2,\ldots,g_n)$ any permutation of $(x_1,x_2,\ldots,x_n)$.
In the case of arbitrary primitive elements
$g_i=q_i(x_1,x_2,\ldots)$ that are homogeneous of a degree $m_i$,
step 2) shows that there are well-defined composition maps for
$\Prim{\mathcal P}$. Thus we have constructed an operad
$\Prim{\mathcal P}$.  The vector space $\Prim(F_{\mathcal
P}(V_X))$ is a $\Prim{\mathcal P}$-subalgebra of $F_{\mathcal
P}(V_X)$ which is free on $V_X$.
 \qed
\end{pf}

\begin{rem}{\rm

The operad $\Prim\As$ is $\Lie$ (cf.\ \cite{lilosci}). This
follows from the Theorem of Friedrichs (cf.\ \cite{lireu}), which
states that Lie polynomials are exactly the polynomials in
non-commuting associative variables which are primitive under
$\Delta_a$. The theorem of Cartier-Milnor-Moore (cf.\
\cite{liqui}, Appendix B), together with the theorem of
Poincar\'e-Birkhoff-Witt, shows that the category of cocommutative
connected graded Hopf algebras is equivalent to the category of
Lie algebras.

The operad $\Dend$ of dendriform algebras of \cite{lilodia} can be
provided with a coherent unit action such that $\vert\prec y=0,\
x\prec \vert=x,\ x\succ \vert=0,\ \vert\succ y=y$, see
\cite{lilosci}.
The free $\Dend$-algebra $(K\YTree^{\infty},
\prec,\succ)$ can be provided with the structure of a $\Dend$-Hopf
algebra (see \cite{lilrhopftree}). It is exactly the $\Dend$-Hopf
algebra structure given by $\Delta_a$. Ronco \cite{liroa} has
determined the primitive elements, the operad $\Prim\Dend$ is the
operad of $\Brace$-algebras, special pre-Lie algebras equipped
with
 $n$-ary operations $\langle \ldots \rangle: A^{\otimes n}\to A$ for each $n$.
For the operad $\Dend$ they play the role of $\Lie$-algebras in
the analogues of Cartier-Milnor-Moore and Poincar\'e-Birkhoff-Witt
theorems (\cite{lirobcr}, \cite{lichamimo}). Further results in
this direction can be found in \cite{lilr04}.

}
\end{rem}

\end{section}

\begin{section}{$\Mgn$-Hopf algebras with duals related to shuffles}

We consider the connected graded $\Mgom$-Hopf algebra
$A=K\{X\}_{\omega}$ with $\Delta_a$, and we consider $K\{X\}_N$ as
a sub $\Mgn$-Hopf algebra for each $N$.

\begin{lem}
The $\Mgn$-Hopf algebras $K\{X\}_N$, $2\leq N \leq \omega$, are
cocommutative (in the sense of Lemma \ref{lemdiagonal}). \qed
\end{lem}

The following Lemma generalizes the formula for $N=2$ given in
\cite{lighmag}.

\begin{lem}
Let $T$ be a tree in $K\{X\}_{\omega}$. For $I$ in ${\rm Le}(T)$
let $I^c$ denote the complement of $I$ in ${\rm Le}(T)$. Then the
image $\Delta_a(T)$ of $T$ is given by the formula
\begin{equation*}
\sum_{I\subseteq {\rm Le}(T)}{\rm red} (T\vert I)\otimes {\rm
red}(T\vert I^c),
\end{equation*}
where the construction of ${\rm red} (T\vert I)$ can be sketched
as follows: All vertices and edges of $T$ that do not lie on some
path from a leaf in $I$ to the root are removed from $T$. Then the
necessary contractions are made to obtain a reduced (admissibly
labeled) tree ${\rm red} (T\vert I)$.
\end{lem}
\begin{pf}
If ${\rm Le}(T)=\emptyset$, the formula says that
$\Delta_a(1)=1\otimes 1$. If $T$ consists of one leaf labeled by
$x_k$, then
\begin{equation*}
\Delta_a(T)=x_k\otimes 1+1\otimes x_k= {\rm red} (T\vert{\rm
Le}(T))\otimes {\rm red}(T\vert\emptyset)+{\rm red}
(T\vert\emptyset)\otimes {\rm red}(T\vert{\rm Le}(T)).
\end{equation*}
Assume now that $T$ has at least two leaves. We may iteratively
apply the $\Mgom$-algebra morphism property
$\Delta_a\circ\vee^n=\vee^n_{A\otimes
A}\circ(\Delta_a\otimes\ldots\otimes\Delta_a)$ of $\Delta_a$ and
expand $\vee^{n_1}_{A\otimes A}(\vee^{n_2}_{A\otimes
A}\ldots\bigl(\ldots, x_k\otimes 1+1\otimes
x_k,\ldots),\ldots\bigr)$ distributively to show that
$\Delta_a(T)$ is given by summands where the tree $T$ is splitted
as indicated, according to the choice of $I\subseteq{\rm Le}(T)$.
\qed
\end{pf}

\begin{defn}
Given two planar admissibly labeled trees $T^1, T^2$ with $n_1,
n_2$ leaves, and a planar tree $T$ with $n_1+n_2$ leaves, we say
that $T$ is a shuffle of $T^1$ and $T^2$ in $K\{X\}_{\omega}$, if
\begin{equation*}
{\rm red}(T\vert I)=T^1, {\rm red}(T\vert I^c)=T^2
\end{equation*}
for some subset $I\subseteq {\rm Le}(T)$.

Similarly defined is the notion of a shuffle in $K\{X\}_N$. If
$T^1, T^2$  both belong to $K\{X\}_N$, $N\in \N$, one can consider
their shuffles in $K\{X\}_N$ and the (larger) set of shuffles in
$K\{X\}_{\omega}$.

 For any tree $T$ in
$K\{X\}_{\omega}$ we define a $K$-linear map $\partial_T:
K\{X\}_{\omega}\to  K\{X\}_{\omega}$ by
\begin{equation*}
\Delta_a(f)=\sum_{T\in\PRTrx} T\otimes
\partial_T(f),
\end{equation*}
where  $f\in K\{X\}_{\omega}$. Especially, for
$T=1_{K\{X\}_{\omega}}=\emp$,  $\partial_{\emp}=\id$.
\end{defn}
\begin{rem}{\rm
The map $\partial_T$ may be called a generalized differential
operator. For any $x_k\in X$, $\partial_{x_k}$ is the unique
mapping
 $D:K\{X\}_{\omega}\to K\{X\}_{\omega}$ with
$D(x_l):=\begin{cases} 1 &:\ k=l \\ 0 &:\ k\neq l\\
\end{cases}$
satisfying the Leibniz rule
\begin{equation*}
D(\vee^n(v_1,\ldots,v_n))=\sum_{i=1}^n
\vee^{n}(v_1,\ldots,v_{i-1},D(v_i),v_{i+1},\ldots,v_n) \text{ for
all } v_i.
\end{equation*}
}
\end{rem}

\begin{prop}
\begin{itemize}
\item[(i)] For $S, T\in\PRTrx$, $T=\vee^p(T^1.T^2\ldots T^p)$, we have
\begin{equation*}
\partial_{S}(T)=
\sum_{\vee^p(S^1.S^2\ldots
S^p)=S}\vee^p\bigl(\partial_{S^1}(T^1)\ldots\partial_{S^p}(T^p)\bigr),
\end{equation*}
where the sum is over all not necessarily non-empty trees $S^1,
S^2,\ldots, S^p$ such that $\vee^p(S^1.S^2\ldots S^p)=S$.
\item[(ii)]
Let $x_k\in X$, and let $W(k,n)$ be the set of reduced trees with
$n$ leaves that are all labeled by $x_k$. Then:
\begin{equation*}
\sum_{S \in W(k,n)}\partial_S=\frac{1}{n!}(\partial_{x_k})^n.
\end{equation*}
\end{itemize}
\end{prop}
\begin{pf}
Since $\Delta_a$ is a $\Mgom$-algebra morphism, we get assertion
(i).
 To prove the equation in assertion (ii), we apply both sides
  to $T=\vee^p(T^1\ldots T^p)$. The assertion is
trivial for $n=0$, $n=1$, or if $T$ has less than $n$ leaves. On
the one hand we have
\begin{equation*}
\frac{1}{n!}\bigl(\partial_{x_k}\bigr)^n(T)=\frac{1}{n!}\sum_{i_1+\ldots+i_p=n}
{n \choose
i_1,\ldots,i_p}\vee^p\bigl((\partial_{x_k})^{i_1}(T^1)\ldots
(\partial_{x_k})^{i_p}(T^p)\bigr).
\end{equation*}
On the other hand, by (i),
\begin{equation*}
\sum_{S \in W(k,n)}\partial_S(T) =\sum_{i_1+\ldots+i_p=n}\ \
\sum_{S^j \in W(k,\
i_j)}\vee^p\bigl(\partial_{S^1}(T^1)\ldots\partial_{S^p}(T^p)\bigr),
\end{equation*}
and assertion (ii) follows by induction. \qed
\end{pf}

\begin{lem}\label{critpseudo}
For $f$ homogeneous of degree $n$, $f$ is primitive if and only if
$\partial_T(f)=0$ for all monomials $T\in K\{X\}_{\omega}$ with
$1\leq \deg T < \frac{n+1}{2}$.
\end{lem}
\begin{pf}
By definition $f$ is primitive if and only if $\partial_T(f)=0$
for all monomials $T\in K\{X\}_{\omega}$.
 Using the cocommutativity of $\Delta_a$, the criterion follows.
\qed
\end{pf}

The commutator operation $[x,y]:=\vee^2(x,y)-\vee^2(y,x)\in
\Prim\Mgom(2)$ is only the first in a large number of primitive
operations, see Example \ref{exprim}.

In view of Lemma\ \ref{lemdualhopf}, we can describe the graded
dual $K\{X\}^{*g}_{\omega}$ of the $\Mgom$-Hopf algebra
$K\{X\}_{\omega}$ equipped with $\Delta_a$ as follows. In analogy
to the classical case, see Remark\ \ref{remcocom}, the commutative
associative binary operation corresponding to $\Delta_a$ is called
planar tree shuffle multiplication (or planar shuffle product),
and the coefficients ${ T \choose T^1, T^2}$ are called planar
binomial coefficients, see \cite{ligepbc}.

\begin{prop}\label{propdualshuffle}

The vector spaces $K\{X\}_{\omega}$ and $K\{X\}^{*g}_{\omega}$ can
be identified by mapping the basis given by trees $T$ on the
corresponding dual basis elements $\delta_T$. Then the commutative
associative multiplication $\Delta_a^{*}$ is given by the binary
operation $\sh: K\{X\}_{\omega}\otimes K\{X\}_{\omega}\to
K\{X\}_{\omega}$ induced by
\begin{equation*}
T^1 \otimes T^2 \mapsto \sum_{\text{all shuffles } T \text{ of }
T^1 \text{ and } T^2} { T \choose T^1, T^2}\ T\text{ \ \ for tree
monomials } T^1, T^2,
\end{equation*}
where ${ T \choose T^1, T^2}\geq 1$ is the number of subsets
$I\subseteq {\rm Le}(T)$ with
\begin{equation*}
{\rm red}(T\vert I)=T^1, {\rm red}(T\vert I^c)=T^2.
\end{equation*}
Especially $\ 1\sh 1=1$. \qed
\end{prop}

\begin{rem}{\rm
\begin{itemize}
\item[(i)]
In the case where $T^1$ has $k$ leaves labeled bijectively by
$x_1$, $\ldots$, $x_k$, and $T^2$ has $n-k$ leaves labeled
bijectively by $x_{k+1},\ldots, x_n$, then no coefficient ${ T
\choose T^1, T^2}>1$ can occur. We get a generalization of the
well-known shuffle multiplication of permutations
$\sh:K\Sigma_k\times K\Sigma_{n-k}\to K\Sigma_n$.
\item[(ii)]
The graded dual of the $\Mgn$-Hopf algebra $K\{X\}_N$ equipped
with $\Delta_a$ is the quotient of $(K\{X\}_{\omega},\sh)$ with
respect to the projection $K\{X\}_{\omega}\to K\{X\}_N$ which is
the identity on the subspace $K\{X\}_N$ and 0 on its complement.
By abuse of notation, we denote all these shuffle multiplications
by $\sh$.
\end{itemize}

}
\end{rem}

\begin{exmp}{\rm
The product $x_1\sh x_2 \sh \ldots \sh x_n$ is given by
\begin{equation*}
\sum_{T\in\PRTree^n}\ \ \sum_{\sigma\in\Sigma_n}\ T^{\sigma}
\end{equation*}
where $T^{\sigma}$ is the admissibly labeled tree with first leaf
labeled by $x_{\sigma(1)}$, second by $x_{\sigma(2)}$, and so on.
This can be shown by induction (the case $n=1$ being trivial)
using the fact that every term of $x_1\sh x_2 \sh \ldots \sh x_n$
occurs in a unique way as a shuffle of a term of $x_1\sh x_2 \sh
\ldots x_{n-1}$ and $x_n$.

}
\end{exmp}

\end{section}

\begin{section}{An analogon of Poincar\'e-Birkhoff-Witt}

As a vector space (in fact as a coalgebra) the free $\As$-algebra
on $V$ is isomorphic to the free $\Com$-algebra generated by all
Lie polynomials, i.e.\ by the primitive elements. This is the
Poincar\'e-Birkhoff-Witt theorem.

To describe the operads $\Prim\Mgn$, we are going to use an
analogon of this theorem.

First, we need to describe the primitive elements as irreducible
elements with respect to the shuffle multiplication, and we also
need a description of the operation $(\vee^2)^*:K\{X\}_{\omega}\to
K\{X\}_{\omega}\otimes K\{X\}_{\omega}$.

\begin{prop}\label{propshuffortho}

Let $f\in K\{X\}_{N}^{(n)}$ be homogeneous of degree $n\geq 1$.
\begin{itemize}
\item[(i)]
If $f$ is primitive, then $\langle f, g_1\sh g_2 \rangle=0$ for
all homogeneous $g_1, g_2\in K\{X\}_{N}$ of degree $\geq 1$.
\item[(ii)]
If $S$ is an admissibly labeled tree of degree $1\leq k \leq n-1$,
then $\partial_S(f)=0$ if and only if $\langle f, S\sh g
\rangle=0$ for all $g\in K\{X\}_{N}^{(n-k)}$.
\item[(iii)]
If the homogeneous element $f$ of degree $n$ is
 orthogonal to all shuffle products $S\sh T$ of trees
$S\in K\{X\}_{N}^{(k)}$, $T\in K\{X\}_{N}^{(n-k)}$  with $1\leq k
< \frac{n+1}{2}$, then $f$ is primitive.
\end{itemize}
\end{prop}

\begin{pf}
We can consider elements of $K\{X\}_{\omega}$. By definition,
$\langle \Delta_a(f),  g_1\otimes g_2 \rangle=\langle f, g_1\sh
g_2 \rangle.$
Clearly, if $f$ is primitive, $\langle f, g_1\sh
g_2\rangle=\langle f,g_1\rangle\langle 1,g_2 \rangle +\langle
1,g_1\rangle\langle f,g_2 \rangle=0$.

If $\partial_S(f)\neq 0$ for some tree monomial $S$ of degree
$1\leq k \leq n-1$, then by definition
\begin{equation*}
0\neq \langle \Delta_a(f),  S\otimes g \rangle=\langle f, S\sh g
\rangle \text{ for some } g\in K\{X\}_{\omega}^{(n-k)}.
\end{equation*}
If $f$ is orthogonal to all shuffle products $S\sh S'$, $S, S'$
trees with $\deg S=k \geq 1$, then
\begin{equation*}
\Delta_a(f)=\sum_{T\in\PRTrx-\{S\}} T\otimes
\partial_T(f)
\end{equation*}
and $\partial_S(f)=0$. Then assertion (iii) follows by Lemma\
\ref{critpseudo}. \qed
\end{pf}

\begin{defn}

Let $N\in \N_{\geq 2}\cup\{\omega\}$. The maps given by
$(\vee^k)^*:K\{X\}_{N}^{*g}\to (K\{X\}_{N}^{*g})^{\otimes k}$ are
simply denoted by $\nabla_k$, for each $k\leq N$, see Lemma\
\ref{lemdualhopf}.

Especially, we consider $\nabla_2$ as a map $K\{X\}_{N}\to
K\{X\}_{N}\otimes K\{X\}_{N}$ for each $N$.

\end{defn}

\begin{lem}\label{lemcomag}

We consider the unitary $\Com$-algebras $A=(K\{X\}_N,\sh)$, $N\in
\N_{\geq 2}\cup\{\omega\}$. By $\bar A$ we denote the augmentation
ideal. Then:
 \begin{itemize}
\item[(i)] The $K$-linear map $\nabla_2:A \to A\otimes A$ is a morphism
of $\Com$-algebras.
\item[(ii)] The $K$-linear map $\nabla_2$ is non-coassociative. If $T$ is
an admissibly labeled tree, then
\begin{equation*}
\nabla_2(T)=\begin{cases} T\otimes 1+1\otimes T +T^1\otimes T^2&:
T=\vee^2(T^1.T^2)\\ T\otimes 1+1\otimes T &: \text{ar}_\rho\neq 2\
\ (\rho \text{ the root}).\\ \end{cases}
\end{equation*}
\item[(iii)] If $\nabla'_2(f):=\nabla_2(f)-f\otimes 1-1\otimes f$, then
$\nabla'_2(f)\in \bar A\otimes \bar A$ for $f\in \bar A$.
\end{itemize}
In other words, these unitary $\Com$-algebras are co-${\bf D}$
objects, where ${\bf D}$ is the category of unitary magmas.
\end{lem}

\begin{pf}
1) The chosen unit actions on $\Mgn$ and $\Com$ have the property
that the operations on the tensor product $A\otimes A$ are defined
component-wise. We have to check that
\begin{equation*}
\nabla_2\circ\sh=\underbrace{(\sh\otimes\sh)\circ\tau_2}_{\sh_{A\otimes
A}}\circ(\nabla_2\otimes\nabla_2).
\end{equation*}
By looking at the graded dual, this is equivalent to the equation
\begin{equation*}
\Delta\circ\vee^2=(\vee^2\otimes\vee^2)\circ\tau_2\circ(\Delta\otimes\Delta).
\end{equation*}
The latter equation is fulfilled, because
$(\vee^2\otimes\vee^2)\circ\tau_2$ is $\vee^2_{A\otimes A}$ and
$\Delta: A\to A\otimes A$ is a morphism of $\Mag$-algebras.
\newline
2) Since
$\nabla_2$ is determined by the equation
\begin{equation*}
\langle \vee^2(f_1.f_2), g\rangle= \langle f_1\otimes f_2,
\nabla_2(g)\rangle
\end{equation*}
we conclude that, for $T$ an admissibly labeled tree,
\begin{equation*}
\nabla_2(T)=T\otimes 1+1\otimes T
+\sum_{T=\vee^2(T^1.T^2)}T^1\otimes T^2.
\end{equation*}
It also follows that $\nabla_2$ is not coassociative. We note
that, for $f, g\in A$,
\begin{equation*}
\langle f, g\rangle= \langle \vee^2(1.f), g\rangle= \langle
1\otimes f,  \nabla_2(g)\rangle= \langle f\otimes 1,
\nabla_2(g)\rangle.
\end{equation*}
Thus assertion (iii) follows.
 \newline
 3) The categorical
coproduct for $\Com$-algebras is the tensor product $\otimes$. By
(i) and (iii) this means that $\nabla_2$ provides the unitary
$\Com$-algebras $A=(K\{X\}_{N},\sh)$ with the structure of a
co-magma object, with counit given by the augmentation map.
\qed
\end{pf}

\begin{rem}{\rm

The vector space $K\{X\}_{N}$ equipped with $\nabla_2$ is a
non-associative coalgebra in the sense of \cite{ligr}.

}
\end{rem}

\begin{thm}\label{thmpbw}
Let,  for $N\in \N_{\geq 2}\cup\{\omega\}$, $A=K\{X\}_N$ be the
unitary free $\Mgn$-algebra. Let $W=\Prim A$ be its space of
primitive elements, graded by $W^{(k)}=\Prim A^{(k)}$.
\begin{itemize}
\item[(i)]
The unitary $\Com$-algebra $(A,\sh)$ is freely generated by
$W=\Prim A$ with respect to the shuffle multiplication.
\item[(ii)]
Let $n\geq 2$, and let $B^{(n)}$ be the orthogonal complement of
$\Prim A^{(n)}$ in the vector space $A^{(n)}$. If $f_1,\ldots,f_r$
is a basis (consisting of homogeneous elements) of
$\oplus_{k=1}^{n-1}\Prim A^{(k)}$, then a basis of $B^{(n)}$ is
given by the elements
\begin{equation*}
f_{i_1}\sh f_{i_2}\sh\ldots\sh f_{i_k}, 1\leq i_1\leq i_2\leq
\ldots \leq i_k\leq r, \text{ with } \sum_{j=1}^k\deg f_{i_j}=n.
\end{equation*}
\end{itemize}
\end{thm}

\begin{pf}
1) We have shown in Proposition\ \ref{propshuffortho} that the
homogeneous primitive elements of $A$ are
 also the homogeneous $\sh$-irreducible elements of $A$,
i.e.\  the elements $f$ not of the form $g_1\sh g_2$ (for $f$,
$g_1, g_2$ homogeneous of degree $\geq 1$). Thus the unitary
$\Com$-algebra morphism $\pi: K1\oplus F_{\Com}(W)\to (A,\sh)$ is
surjective, and for every proper subspace $U$ of $W$, $K1\oplus
F_{\Com}(U)$ cannot be isomorphic to $(A,\sh)$.
\newline
2) We have to show that $\pi$ is injective, i.e.\ that
 $(A,\sh)$ is a free $\Com$-algebra.
\newline
By Lemma\ \ref{lemcomag}, there exists a morphism $\nabla_2:A \to
A\otimes A$ of $\Com$-algebras, which provides the unitary
$\Com$-algebras $(A,\sh)$ with the structure of a co-${\bf D}$
object in the category of $\Com$-algebras, where ${\bf D}$ is the
category of unitary magmas.
\newline
Over a field $K$ of characteristic 0, all connected (i.e.\
$A^{(0)}=K$) $\Com$-algebras that are equipped with the structure
of a unital co-magma are free. This is the Leray theorem, see
\cite{lioud}. Thus $\pi$ is an isomorphism.
\newline
3)
 Since the
space $\Prim A^{(n)}$ is orthogonal (with respect to $\langle,
\rangle$) to the shuffle products in $A^{(n)}$, see Proposition\
\ref{propshuffortho}, assertion (ii) follows from assertion (i).
\qed
\end{pf}

Dualizing the statement of Theorem\ \ref{thmpbw}(i) and its proof,
we obtain:

\begin{thm}\label{thmpbwdual}
The $\Mgn$-Hopf algebras given by\ \  the free $\Mgn$-algebras
equipped with $\Delta_a$  are cofree co-nilpotent
$\Com$-coalgebras, co-generated by their primitive elements. \qed
\end{thm}

\begin{rem}\label{remusualhopf}{\rm
\begin{itemize}
\item[(i)]
Given an arbitrary graded connected $\Mgn$-Hopf algebra $A$, one
can consider it as a $\Prim\Mgn$-algebra, and its primitive
elements as a ($\Prim\Mag$-)subalgebra. For $K$ of characteristic
0, the situation is analogous to the classical
Cartier-Milnor-Moore theorem (cf.\ \cite{liqui}, Appendix B), and
the Hopf algebra $A$ is of the form $U(\Prim A)$.
 Of course the concrete description of the
functor $U$ depends on a description of the
$\Prim\Mag$-operations.
\newline
Theorem\ \ref{thmpbw}(ii) shows that one may search for orthogonal
projectors
\begin{equation*}
e_n^{(1)},e_n^{(2)},\ldots,e_n^{(n)}: K\{X\}_{N}^{(n)}\to
K\{X\}_{N}^{(n)}
\end{equation*}
which are similar to the Eulerian idempotents (cf.\
\cite{lilohaus}). The idempotent $e_n^{(i)}$ projects elements of
$K\{X\}_{N}^{(n)}$ into the subspace of $i$-factor shuffle
products of primitive elements. Here any series like the
exponential series of \cite{ligeplaexp} which maps primitive
elements on group-like elements is useful.
\item[(ii)]
Theorem\ \ref{thmpbwdual} suggests the definition of a family of
cocommutative $\As$-Hopf algebras $H_N$ related to the given
$\Mgn$-Hopf algebras. The idea is to have a free Lie algebra on a
space $C$ of graded generators associated to $H_N$ by the
 Cartier-Milnor-Moore theorem, such that the
cocommutative coalgebra structure of $H_N$ is the same as the one
on the corresponding $\Mgn$-Hopf algebra.\newline
 One might consider
the $\Prim\Mgn$-operations $[\ldots[x_n, x_{\delta 1}], x_{\delta
2}]\ldots x_{\delta(n-1)}]$, $\delta\in\Sigma_{n-1}, n\geq 2$, and
choose $C$ such that these operations evaluated on $C$ yield all
primitive elements.
\end{itemize}
}
\end{rem}
\end{section}

\begin{section}{The generating series and
representations}

The generating series $f^{\mathcal P}(t)$ of an
 operad ${\mathcal P}$ is the series
$\sum_{n\geq 1}\frac{\dim {\mathcal P}(n)}{n!}t^n$.

 The classical
Poincar\'e-Birkhoff-Witt theorem implies that the generating
series of $\Prim\As=\Lie$ is $f^{\Lie}(t)=\sum_{n\geq
1}\frac{(n-1)!}{n!}t^n= -\log(1-t)$, because its composition with
$f^{\Com}=\exp(t)-1$ is $f^{\As}(t)=\frac{1}{1-t}-1$.

Similarly, it is implied by Theorem\ \ref{thmpbw} that the
generating series of the operads $\Prim\Mgn$ are the logarithms
$\log(1+t)$ of the generating series of the corresponding $\Mgn$.
In fact this holds on the level of characteristic functions in the
ring $\Lambda$ of symmetric functions, equipped with the plethysm
$\circ$. For the theory of symmetric functions, see \cite{limd}.
 Let $p_r, r\geq 1$ denote
the $\Q$-basis of $\Lambda$ given by the power sum symmetric
functions $\sum_{i\geq 1} x_i^r$.

\begin{cor}\label{corlogcatalan}
\begin{itemize}
\item[(i)]
The characteristic of the $\Sigma_n$-module $\Prim\Mag(n)$ is
given by
\begin{equation*}
{\rm ch}_n\bigl(\Prim\Mag(n)\bigr)= \frac{1}{n}\sum_{d\vert
n}\mu(d) c'_{\frac{n}{d}}p_d^{\frac{n}{d}}
\end{equation*}
and the characteristic function ${\rm ch}\bigl(\mathcal P\bigr)=
\sum_{n\geq 1}{\rm ch}_n\bigl({\mathcal P}(n)\bigr)$ for
$\Prim\Mag$ is given by
\begin{equation*}
\sum_{d\geq 1}\frac{\mu(d)}{d} \sum_{k\geq 1} \frac{c'_k}{k}p_d^k.
\end{equation*}
Here $c'_k$ is the $k$-th $\log$-Catalan number.
\item[(ii)]
Especially it holds that the dimension of ${\Prim\Mag}(n)$ is
given by
\begin{equation*}
\dim{\Prim\Mag}(n)=(n-1)! c'_n
\end{equation*}
and the generating series $f^{\Prim\Mag}(t)$ is given by
\begin{equation*}
\log\bigl(\frac{3-\sqrt{1-4t}}{2}\bigr).
\end{equation*}
\item[(iii)]
The analogous assertions of (i) and (ii) hold for the operads
$\Prim\Mgn$, $N\in \N_{\geq 2}\cup\{\omega\}$: The sequence $c'_n$
has to be replaced by the sequence $c[N]'_n$. Especially, for
$\Mgom$ the logarithmic derivative $C'_n$ of the super-Catalan
numbers $C_n$ has to be taken.
 The generating series $f^{\Prim\Mgom}(t)$ is given by
\begin{equation*}
\log\Bigl(\frac{5+t-\sqrt{1-6t+t^2}}{4}\Bigr).
\end{equation*}
\end{itemize}
\end{cor}

\begin{pf}
We may modify a computation given in \cite{liget}. While the
classical Poincar\'e-Birkhoff-Witt theorem implies that ${\rm
ch}(\Lie)={\rm Log}\bigl(1+{\rm ch}(\As)\bigr)$, Theorem\
\ref{thmpbw} implies that ${\rm ch}(\Prim\Mag)={\rm
Log}\bigl(1+{\rm ch}(\Mag)\bigr).$ Here the operation ${\rm Log}$
on symmetric functions is given by
\begin{equation*}
{\rm Log}(1+f)=\sum_{d\geq 1}\frac{\mu(d)}{d}\sum_{n\geq
1}\frac{(-1)^{n+1}}{n}p_d^n\circ f.
\end{equation*}
Clearly ${\rm ch}(\Mag)$ is given by $\sum_{k\geq 1}c_k p_1^k$,
thus $\sum_{n\geq 1}\frac{(-1)^{n+1}}{n}p_d^n\circ {\rm ch}(\Mag)$
is equal to $\log(1+p_d\circ \sum_{k\geq 1}c_k
p_1^k)=\log(1+\sum_{k\geq 1}c_k p_d^k).$
\newline
By definition of the log-Catalan number $c'_k$ as the coefficient
of $t^{k-1}$ in the logarithmic derivative of $\sum_{k\geq 1}c_k
t^k$, we get that $\log(1+\sum_{k\geq 1}c_k t^k)=\sum_{k\geq 1}
\frac{c'_k}{k}t^k.$ Hence we get the asserted expressions for
${\rm ch}(\Prim\Mag)$ and ${\rm ch}_n(\Prim\Mag)$ in (i).

To show assertion (ii), one may repeat the same computation for
generating functions instead of characteristic functions, or apply
the rank morphism that maps $p_1\mapsto t$ and $p_n\mapsto 0,
n>1$.

The same arguments, with Catalan numbers $c_n$ replaced by the
sequence $c[N]_n$, apply to $\Prim\Mgn$, $N\in \N_{\geq
2}\cup\{\omega\}$. \qed
\end{pf}

\begin{rem}{\rm

The occurrence of log-Catalan numbers in dimension formulas for
homogeneous elements of $\Prim K\{X\}$ has independently also been
observed by Bremner, Hentzel, and Peresi in \cite{libhp}. They
consider a set $X$ of $r_i$ generators of multi-degree $e_i$,
$i=1,\ldots,N$, and show that the Witt dimension formula
\begin{equation*}
\frac{1}{\vert{\bf n}\vert}\sum_{k,{\bf d} \atop k {\bf d}={\bf
n}} \mu(k)c'_{\vert{\bf d}\vert}{\vert{\bf d}\vert\choose
d_1,\ldots,d_N}r_1^{d_1}\cdots r_N^{d_N}
\end{equation*}
yields the dimension of the space of homogeneous primitive
elements of multi-degree ${\bf n}= (n_1,n_2,\ldots,n_N)$.

}
\end{rem}

\begin{rem}{\rm

Using Corollary\ \ref{corlogcatalan} it is easy for small $n$ to
describe
 the representation of $\Sigma_n$ given by $\Prim\Mgn(n)$ in terms
 of irreducible representations. In the
basis of Schur functions, one checks that ${\rm ch}_3$ is
$s_{3}+3s_{2,1}+s_{1,1,1}$ for $\Prim\Mag$, and
$2s_{3}+5s_{2,1}+2s_{1,1,1}$ for $\Prim\Mgom$. We get
$3s_{4}+10s_{3,1}+6s_{2,2}+10s_{2,1,1}+3s_{1,1,1,1}$ for
$\Prim\Mag(4)$, and
$8s_{4}+25s_{3,1}+16s_{2,2}+25s_{2,1,1}+8s_{1,1,1,1}$ for
$\Prim\Mgom(4)$.

These representations of $\Sigma_n$ are not given by copies of the
$(n-1)!$-dimensional representation $\Lie(n)$. But they occur as
the representations of primitives associated to the cocommutative
Hopf algebras mentioned in Remark\ \ref{remusualhopf}(ii), because
the homogeneous components of degree $n$ of these Hopf algebras
are $c[N]_n$ copies of the regular representation (in the
$\Mgn$-case). Thus one may obtain this type of representations,
starting with free $\As$-algebra generators, as was pointed out to
us by J.-C.\ Novelli. Here one needs $e[N]_n$ generators in degree
$n$, where
  the generating series $e(t)$ and $c(t)$ for $e[N]_n$ and $c[N]_n$ are
related by $c(t)=\frac{1}{1-e(t)}-1$.

A comment made by F.\ Chapoton is that it would be nice to have
presentations of the operads $\Prim\Mgn$ as Hadamard products, and
one may ask for a possible (anti-)cyclic operad structure, see
\cite{lichacyc}.

An intrinsic characterization of $\Prim\Mgn$ by generators and
relations would naturally be interesting, a question that was
posed (for primitive elements in free $\Mag$-algebras) by Umirbaev
and Shestakov in \cite{lisu}, see also \cite{lighmag}.

}
\end{rem}

\begin{exmp}\label{exprim}{\rm
The commutator operation $\vee^2(x,y)-\vee^2(y,x)\in
\Prim\Mag(2)$\ $=\Prim\Mgom(2)$ is denoted by $[x,y]$. Let
$\langle x,y,z\rangle := (x,z,y)-(x,y,z)$, where $(x,y,z)$ is the
associator operation $\vee^2(\vee^2(x,y),z) -
\vee^2(x,\vee^2(y,z))$. Furthermore, let $\{ x,y,z\}:= (x,z,y)+
(x,y,z).$

The space $\Prim\Mag(3)$ has dimension 8 and is generated (as a
$\Sigma_3$-module) by the operations $\langle x,y,z\rangle$, $\{
x,y,z\}$, and $[[x,y],z]$, which fulfill
\begin{equation*}
\begin{split}
&\langle x,z,y\rangle = -\langle x,y,z\rangle,\  \{ x,z,y\}=\{
x,y,z\},\  [[y,x],z]=-[[x,y],z], \text{and}\\
&\text{ the equation  } \sum_{a,b,c \text{ cyclic}}\langle
a,b,c\rangle= \sum_{a,b,c \text{ cyclic}}[[a,b],c],
\\
\end{split}
\end{equation*}
which is called the non-associative Jacobi relation (cf.\
\cite{lilosci}) or Akivis relation (cf.\ \cite{lisu}). Here the
operation $ \langle x,y,z\rangle+[[y,z],x]$ generates only a
2-dimensional $\Sigma_3$-module.

The space $\Prim\Mag_3(3)=\ldots=\Prim\Mgom(3)$ has dimension 14
and is generated by the operations above together with the
operations $(x,y,z)_t:=\vee^2(\vee^2(x,y),z)-\vee^3(x,y,z)$.

The space $\Prim\Mag(4)$ of dimension of  $3!\cdot c'_3=78$ can be
generated by the operations $p$ and $q$ given by $p(x,t,y,z):=(
xt,y,z) -x\cdot (t,y,z) -t\cdot(x,y,z)$ and $q(x,t,y,z):=( x, ty,
z) -y\cdot (x,t,z)-t\cdot(x,y,z)$.

}
\end{exmp}

\begin{defn}(cf.\ \cite{lisu}, \cite{liper})

A vector space $V$ (over a field $K$ of characteristic 0) together
with $(m+2)$-ary operations $\langle \ {\bf x}\ \vert\ y\ \vert\
z\rangle$, all $m\geq 0$, and $(m+n)$-ary operations $\Phi({\bf
x}\ \vert\ {\bf y})$, all $m\geq 1, n\geq 2$, is called a Sabinin
algebra, if relations hold which can be abbreviated -- using
Sweedler's notation for $\Delta_a(x_1.x_2\ldots x_r)$ -- by:
\begin{equation*}
\begin{split}
&\langle {\bf x}\vert a \vert b\rangle=-\langle {\bf x}\vert
b\vert a\rangle,\\
& \Phi({\bf x}\vert {\bf y})=\Phi(x_{\tau 1},\ldots,x_{\tau
m}\vert y_{\delta
1},\ldots,y_{\delta n}), \text{ all }\tau\in\Sigma_m,\delta\in\Sigma_n,\\
\\
&\langle x_1,\ldots x_r, a,b, z_{r+1},\ldots z_m\vert c \vert
d\rangle-\langle x_1,\ldots x_r, b,a, z_{r+1},\ldots z_m\vert c
\vert d\rangle\\&+ \sum \langle {\bf x}_{(1)}\ \langle {\bf
x}_{(2)}\vert a\vert b\rangle\ z_{r+1}\ldots z_m\ \vert\ c \vert
d\rangle=0, \text{ all } r,
\\
&\sum_{a,b,c \text{ cyclic}}\bigl(\langle {\bf x}c\ \vert\ a\
\vert\ b\rangle+\sum \langle {\bf x}_{(1)}\ \vert\ \langle {\bf
x}_{(2)}\ \vert\ a\ \vert \ b\rangle \ \vert\ c\rangle\bigl)=0.
\\
\end{split}
\end{equation*}
\end{defn}

\begin{rem}{\rm

We used $p$ and $q$ to generate $\Prim\Mag(4)$. More generally,
 Shestakov and Umirbaev \cite{lisu} show that the recursively
 defined operations
\begin{equation*}
P(x_1,\ldots,x_m; y_1,\ldots,y_n; z)=({\bf x},{\bf y},z) -{\sum}'
{\bf x}_{(1)}\ {\bf y}_{(1)}\  P\bigl({\bf x}_{(2)}; {\bf
y}_{(2)};z\bigr),
\end{equation*}
(together with $[x,y]$) give a complete set of primitive
operations. Here products of more than two arguments are evaluated
from left to right, and $\sum'$ indicates that there are no terms
with ${\bf x}_{(i)}={\bf y}_{(i)}=1$. Moreover they show that the
space $B$ of primitive elements together with
\begin{equation*}
\begin{split}
&\langle y \vert z\rangle := -[y,z],\ \ \text{and for } m\geq 1:\
\langle {\bf x} \vert y \vert z\rangle :=
P({\bf x};z;y)-P({\bf x};y;z),\\
&\text{ and for } m\geq 1, n\geq 2:\ \Phi({\bf x}\vert {\bf y}) :=
 \sum_{\tau,\delta\in\Sigma}\frac{P(x_{\tau 1},\ldots,x_{\tau m}; y_{\delta 1},\ldots,y_{\delta
(n-1)}; y_{\delta n})}{m! n!}
\\
\end{split}
\end{equation*}
is a Sabinin algebra $G(B)$. Especially this means that these
(multilinear) Sabinin operations form a sub-operad $\Sab$ of
$\Prim\Mag$.
\newline
We denote  $m! n!\Phi({\bf x}\vert {\bf y})$ also by $\{ {\bf x} \
\vert\  {\bf y}\}$.
 }
\end{rem}

\begin{conj}
The operads $\Sab$ and
 $\Prim\Mag$ are the same.

\end{conj}

This conjecture is enforced by recent results of P\'erez-Izquierdo
\cite{liper}, which show that Sabinin algebras $V$ have a
universal enveloping ($\Mag$-)algebra with a
Poincar\'e-Birk\-hoff-Witt basis (given by monomials on $V$) and
that there is a Cartier-Milnor-Moore theorem for Sabinin algebras.
The free $\Mag$-algebra also has a Poincar\'e-Birk\-hoff-Witt
basis formed by primitive elements (by \cite{lisu} or Theorem
\ref{thmpbwdual}).

We know that $\Prim\Mag(n)=\Sab(n)$ for $n=2,3.$ In the degree 4
component of $\Prim\Mag$ of dimension 78 the following Sabinin
relation holds:
\begin{equation*}
\sum_{a,b,c \text{ cyclic}}\langle x,c \ \vert\  a\ \vert\
b\rangle = \sum_{a,b,c \text{ cyclic}}\bigl( \langle
 x,[a,b],c\rangle+[\langle x,a,b\rangle,c]\bigr),
\end{equation*}
where $\langle x , t\ \vert\ y\ \vert\ z\rangle :=
p(x,t,z,y)-p(x,t,y,z)$. By a computation we find that the
$\Sigma_4$-subspace given by the iterated commutators and the
operations $[\langle x,t,y\rangle,z]$, $[\{x,t,y\},z]$, $\langle
[x,t],y,z\rangle$, $\langle x,[t,y],z\rangle,$ $\{ [x,t],y,z\}$,
and $\{ x,[t,y],z\}$ has dimension 65. Three more basis elements
of $\Prim\Mag(4)$ may be given by
 $\langle x , t\
\vert\ y\ \vert\ z\rangle$,
 $\langle y , z\
\vert\ x\ \vert\ t\rangle$, and
 $\langle z , t\
\vert\ y\ \vert\ x\rangle$.
 Together with the (independent) 6+4 operations
given by
 $\{x t \vert y z\}$ and  $\{x \vert tyz\}$, we arrive at a description of
  $\Prim\Mag(4)$ by Sabinin operations, confirming that
$\Prim\Mag(4)=\Sab(4).$

\begin{rem}{\rm

The commutative version of the operad $\Mag$ is the operad $\Cmg$
of commutative magma algebras. It is the free operad generated by
a free binary commutative operation. The trees symbolizing
operations are abstract instead of planar trees. If we replace
planar trees by abstract trees, we get 'commutative versions'
$\Cmgn$ of the operads $\Mgn$. In this paper we have focused on
the operads $\Prim\Mgn$. The same techniques can be applied to
describe operads $\Prim\Cmgn$, and we will consider these operads
in a subsequent paper. For $\mathcal P=\Prim\Cmg$, we get that
${\mathcal P}(2)=0$ and ${\mathcal P}(3)$ is 2-dimensional,
generated by the associator operation $(x_1,x_2,x_3)=(x_1\cdot
x_2)\cdot x_3 - x_1\cdot (x_2\cdot x_3) =-(x_3,x_2,x_1)$ subject
to the relation $(x_1,x_2,x_3)+(x_2,x_3,x_1)+(x_3,x_1,x_2)=0$. In
terms of irreducible representations (or Schur functions), ch$_3$
is $s_{2,1}$ and ch$_4$ is $s_4+s_{3,1}+s_{2,2}$ (with dimension
6).

}
\end{rem}

\end{section}

\goodbreak

\end{document}